\documentclass[12pt]{article}

\newtheorem{thm}{Theorem}[section]
\newtheorem{prop}[thm]{Proposition}
\newtheorem{lem}[thm]{Lemma}
\newtheorem{dfn}[thm]{Definition}

\newtheorem{con}[thm]{Conjecture}

\newcommand{\qed}{\hfill \fbox{}\medskip}
\newcommand{\proof}{\medskip\noindent{\bf Proof.}\quad }

\title{A Commutative Family of Integral Transformations
and Basic Hypergeometric Series. I. 
Eigenfunctions}

\author{Jun'ichi Shiraishi\\
\\
{\it Graduate School of 
Mathematical Science, }\\ {\it University of Tokyo, Tokyo, Japan}}

\date{}

\begin{document}

\maketitle

\maketitle

\begin{abstract}
It is conjectured that 
a class of $n$-fold integral transformations 
$\{I(\alpha)|\alpha\in {\bf C} \}$
forms a mutually commutative family, namely,  we have
$I(\alpha) I(\beta)=I(\beta) I(\alpha) $ for 
${}^\forall \alpha,{}^\forall \beta \in {\bf C}$.
The commutativity of $I(\alpha)$ 
for the two-fold integral case is proved
by using several summation and transformation formulas
for the basic hypergeometric series. 
An explicit formula for the complete system of the eigenfunctions for 
$n=3$ is conjectured. In this formula and 
in a partial result for $n=4$, it is observed that
all the eigenfunctions do not depend on the spectral parameter $\alpha$
of $I(\alpha)$.

\end{abstract}

\section{Introduction}
It was pointed out in \cite{S} that a certain class of $n$-fold integral
transformations plays an essential role for a 
study of the vertex operator $\Phi(\zeta)$ for Baxter's eight-vertex model \cite{B1,B2,B}.
(As for the definition of the intertwiner $\Phi(\zeta)$, 
see \cite{FIJKMY1,FIJKMY2} and
\cite{S}.)
More precisely, in \cite{S} the matrix elements
$\langle \Phi(\zeta_1)\Phi(\zeta_2)\cdots \Phi(\zeta_n)\rangle$
were represented by applying the $n$-fold integral transformation
to a basic hypergeometric series. 
(See Eq. (47) of \cite{S}.)
The aim of the present paper is to investigate the structure of the integral
transformation to obtain a better understanding of 
the eight-vertex vertex operator $\Phi(\zeta)$.
\bigskip

Let us first recall the notations for the integral representation given in \cite{S}.
Let $h(\zeta)$ and $g(\zeta)$ be the functions defined by
\begin{eqnarray}
h(\zeta)&=&(1-\zeta)
{(q t^{-1} \zeta;q)_\infty \over
(t\zeta;q)_\infty },\\
g(\zeta)&=&
{(q^{1\over 2}t^{1\over 2} \zeta;q)_\infty \over
(q^{1\over 2}t^{-{1\over 2}} \zeta;q)_\infty }.
\end{eqnarray}
Here we have used the standard notation for the
$q$-shifted factorial (\ref{po-inf}).
Note that we have slightly modified the definition of $h(\zeta)$
from the one given by Eq.(4) in \cite{S}. One finds this change in $h(\zeta)$ will 
simplify our discussion given in what follows.
\medskip

Let us introduce a space of formal power series of degree zero in
variables $\zeta_i$ ($i=1,2,\cdots,n$)
corresponding to the positive cone of the $A_{n-1}$ type root lattice
\begin{eqnarray}
{\cal F}_n=
\left\{
\sum_{i_1,i_2,\cdots,i_{n-1}\geq 0}^\infty
c_{i_1,i_2,\cdots,i_{n-1}}\left(\zeta_2\over \zeta_1\right)^{i_1}
\left(\zeta_3\over \zeta_2\right)^{i_2}
\cdots \left(\zeta_n\over \zeta_{n-1}\right)^{i_{n-1}} \right\}.
\end{eqnarray}
Note that the matrix elements of the eight-vertex 
vertex operators belong to this space, namely
$\langle \Phi(\zeta_1)\Phi(\zeta_2)\cdots \Phi(\zeta_n)\rangle \in{\cal F}_n$.
Our task in this paper is to propose an integral operator which acts on ${\cal F}_n$,
and study some basic properties of it.
Applications will be considered 
in the continuations of the present paper \cite{S2,S3}.
\medskip

The central object in the present article is given as an integral
transformation acting on  ${\cal F}_n$.
\begin{dfn}
Let $(s_1,s_2,\cdots,s_n)\in{\bf C}^n$ and $\alpha\in{\bf C}$ be parameters. 
Let $q,t\in{\bf C}$ be parameters satisfying the condition $|qt^{-1}|<1$.
We will assume that all these parameters are generic unless otherwise stated.
The $n$-fold integral transformation
$I(\alpha)=I(\alpha;s_1,s_2,\cdots s_n,q,t)$ is defined by
\begin{eqnarray}
&&I(\alpha)  f(\zeta_1,\cdots,\zeta_n) \nonumber\\
&=&
\prod_{i=1}^n \left( 
{ (q t^{-1} ;q)_\infty \over 
(\alpha s_i^{-1} q t^{-1} ;q)_\infty}
{ (q;q)_\infty \over 
(\alpha^{-1}s_i q;q)_\infty}
\right)\\
&&
\times\prod_{i<j} h(\zeta_j/\zeta_i)
\oint_{C_1}\cdots\oint_{C_n} {d\xi_1 \over 2\pi i\xi_1}\cdots
{d\xi_n \over 2\pi i\xi_n}
\prod_{i=1}^n
{\Theta_{q}(\alpha s_i^{-1} q^{1\over 2}t^{-{1\over 2}} \zeta_i/\xi_i)\over 
\Theta_{q}(q^{1\over 2}t^{-{1\over 2}}\zeta_i/\xi_i)} \nonumber\\
&&
\times
\prod_{k=1}^n
\left[
\prod_{i<k} g(\zeta_k/\xi_i)
\prod_{j\geq k} g(\xi_j/\zeta_k)
\right]
f(\xi_1,\cdots,\xi_n), \nonumber
\end{eqnarray}
where the integration contours $C_i$ are given by the conditions $|\zeta_i/\xi_i|=1$,
and the theta function $\Theta_{q}(\zeta)$ is 
given by (\ref{theta-def}).
\end{dfn}
Then our main statement in this article is 
\begin{con}\label{Con1}
For ${}^\forall \alpha,{}^\forall \beta \in{\bf C}$,
the integral transformations 
$I(\alpha), I(\beta)$ acting on the space ${\cal F}_n$ are mutually commutative
\begin{eqnarray}
I(\alpha) I(\beta)=I(\beta) I(\alpha). \label{II=II}
\end{eqnarray}
\end{con}
Here the parameters $(s_1,s_2,\cdots,s_n)$, $q$ and $t$ should be fixed, namely
Eq.(\ref{II=II}) means 
\begin{eqnarray}
&&I(\alpha;s_1,\cdots,s_n,q,t) 
I(\beta;s_1,\cdots,s_n,q,t)\nonumber\\
&=&
I(\beta;s_1,\cdots,s_n,q,t)
I(\alpha;s_1,\cdots,s_n,q,t).\nonumber
\end{eqnarray}

In Section 2, we will prove the commutativity of the
integral transformation $I(\alpha)$ for the case of $n=2$
by using some summation and transformation formulas for the
basic hypergeometric series.
\begin{thm}\label{Thm1}
For $n=2$, we have
\begin{eqnarray}
I(\alpha) I(\beta)=I(\beta) I(\alpha) ,
\end{eqnarray}
on ${\cal F}_2$.
\end{thm}
In Section 3, the complete system of the eigenfunctions for $n=3$ will be conjectured.
This gives us a strong support for Conjecture \ref{Con1} for $n=3$, since
no dependence on the spectral parameter $\alpha$ is observed in all the eigenfunctions.
Some evidence for the commutativity for $n=4$ will also be given.
\medskip

Conjecture \ref{Con1} means that there exists
a quantum mechanical integrable system whose Hamiltonian
is given by the integral operator $I(\alpha)$.
It is an interesting problem to relate this quantum mechanical 
system with already known one.
This will be considered in the next paper \cite{S2}.
A Macdonald-type difference operator will be introduced and 
its commutativity with the action of $I(\alpha)$ will be discussed.
\medskip

The meaning of the commutativity $[I(\alpha),I(\beta)]=0$, however, 
still remains  unclear from the lattice model point of view.
Firstly, 
we do not understand the role of  Eq.(\ref{II=II})
for the characterization of the correlation functions of the eight-vertex model,
since the integral operator $I(\alpha)$
was heuristically introduced through the investigation based on the
free field representation of $\Phi(\zeta)$.
Next, we lack an explanation for 
Eq.(\ref{II=II}) based on the commutative family 
generated by the row-to row transfer matrix
of the eight-vertex model.
\bigskip

The plan of the paper is as follows. In Section 2, a proof of 
Theorem \ref{Thm1} is given. Several summation and transformation
formulas for the basic hypergeometric series will be used there.
In Section 3, basic properties of the eigenfunctions of $I(\alpha)$
are described. Explicit formulas for the eigenfunctions for $n=3,4$ will be conjectured
and shown that these are independent of the spectral parameter $\alpha$.
For $n=3$, it gives us the complete system of the eigenfunctions.
For $n=4$, however, only several lower terms will be obtained.
Concluding remarks are given in Section 4.
\bigskip

Throughout the paper, we use the standard notations
for the $q$-shifted factorials
\begin{eqnarray}
(a;q)_n&=&
\left\{
\begin{array}{ll}
1 & (n=0), \\[3mm]
(1-a)(1-a q)\cdots(1-a q^{n-1}) &(n=1,2,\cdots),
\end{array}
\right. \\[4mm]
(a_1,a_2,\cdots,a_m;q)_n&=&(a_1;q)_n(a_2;q)_n\cdots (a_m;q)_n,\\[4mm]
(a;q)_\infty &=& \prod_{k=0}^\infty (1-a q^k), \label{po-inf}\\[4mm]
(a;q)_{-n}&=&
{1\over (1-a q^{-1})(1-a q^{-2})\cdots (1-a q^{-n})}
\quad (n=1,2,\cdots),\nonumber\\
\end{eqnarray}
and the basic hypergeometric series
\begin{eqnarray}
{~}_{r+1}\phi_r\left(
{a_1,a_2,\cdots,a_{r+1}\atop
b_1,b_2,\cdots,b_r};q,z
\right)
=
\sum_{n=0}^\infty
{(a_1,a_2,\cdots,a_{r+1};q)_n\over
(b_1,b_2,\cdots,b_{r},q;q)_n}z^n.
\end{eqnarray}
We also use the notation for the elliptic theta function as
\begin{eqnarray}
\Theta_q(z)=(z;q)_\infty(q/z;q)_\infty(q;q)_\infty.\label{theta-def}
\end{eqnarray}
\bigskip

\section{Proof of Theorem \ref{Thm1}}
In this section, we prove Theorem \ref{Thm1}.
For the two-fold integral case, the integral transformation reads
\begin{eqnarray}
&&I(\alpha)  f(\zeta_1,\zeta_2)\nonumber \\
&=&
\prod_{i=1}^2
\left( { (q t^{-1};q)_\infty \over 
(\alpha s_i^{-1}q t^{-1} ;q)_\infty}
{ (q;q)_\infty \over 
(\alpha^{-1}s_i q;q)_\infty}
\right)h(\zeta_2/\zeta_1)\\
&&
\times
\oint \oint {d\xi_1 \over 2\pi i\xi_1} 
{d\xi_2 \over 2\pi i\xi_2}
{\Theta_{q}(\alpha s_1^{-1} q^{1\over 2}t^{-{1\over 2}}\zeta_1/\xi_1)\over 
\Theta_{q}(q^{1\over 2}t^{-{1\over 2}}\zeta_i/\xi_i)}
{\Theta_{q}(\alpha s_2^{-1} q^{1\over 2}t^{-{1\over 2}}\zeta_2/\xi_2)\over 
\Theta_{q}(q^{1\over 2}t^{-{1\over 2}}\zeta_i/\xi_i)}
 \nonumber\\
&&
\times
g(\xi_1/\zeta_1)g(\xi_2/\zeta_1)
g(\zeta_2/\xi_1)g(\xi_2/\zeta_2)
f(\xi_1,\xi_2). \nonumber
\end{eqnarray}
Let us prove the following theorem 
which is equivalent to Theorem \ref{Thm1}.
\begin{thm}\label{Thm2}
For the two-fold integral case, the operator $I(\alpha)$ has 
the eigenvalues and eigenvectors:
\begin{eqnarray}
I(\alpha) f_j(\zeta_1,\zeta_2) =\lambda_j (\alpha) f_j(\zeta_1,\zeta_2),
\end{eqnarray}
for $j=0,1,2,\cdots$, where
\begin{eqnarray}
\lambda_j(\alpha)
&=&
{(\alpha s_1^{-1};q)_{-j} \over 
  (\alpha s_1^{-1} q t^{-1};q)_{-j}}
{(\alpha s_2^{-1};q)_j \over 
  (\alpha s_2^{-1}q t^{-1};q)_j}\\
&=&
{(\alpha^{-1} s_1 t ;q)_j \over (\alpha^{-1} s_1q;q)_j}
{(\alpha s_2^{-1};q)_j \over 
 (\alpha s_2^{-1}q t^{-1};q)_j}
\left(q t^{-1} \right)^j, \label{lambdaj}
\end{eqnarray}
and
\begin{eqnarray}
&&f_j(\zeta_1,\zeta_2)\nonumber\\
&=&
\zeta^j\times
{}_4\phi_3\left(
{  s q^{2j}t^{-1} , s^{1\over 2} q^{j+1} t^{-{1\over 2}},
- s^{1\over 2}  q^{j+1} t^{-{1\over 2}} ,t^{-1}
\atop
s^{1\over 2}  q^{j} t^{-{1\over 2}} ,-s^{1\over 2}  q^{j}t^{-{1\over 2}},
 s q^{2j+1}}
;q, t\zeta
\right)\label{f_j}\\
 &=&
 \zeta^j(1-\zeta)\times
{}_2\phi_1\left(
{  q t^{-1}, s q^{2j+1} t^{-1},
\atop
s q^{2j+1} }
;q, t\zeta
\right).
 \end{eqnarray}
 Here we have denoted $\zeta=\zeta_2/\zeta_1$ and
 $s=s_1/s_2$ for short.
\end{thm}
  
Since the eigenfunctions $f_j(\zeta_1,\zeta_2)$ do not depend on 
the parameter $\alpha$ and
our space of series ${\cal F}_2$ is spanned by $\{f_j(\zeta_1,\zeta_2)|j=0,1,2\cdots\}$,
we immediately see that Theorem \ref{Thm2} is equivalent to 
Theorem \ref{Thm1}.
\medskip

 Note that we have expressed the eigenfunctions $f_j(\zeta_1,\zeta_2)$ in two ways
by using 
\begin{lem}We have
\begin{eqnarray}
(1-z){}_2\phi_1\left(
{a,b\atop aqb^{-1}};q,zqb^{-1}
\right)
=
{}_4\phi_3\left(
{a^{1\over 2}q^{1\over 2},
-a^{1\over 2}q^{1\over 2},aq^{-1},bq^{-1}\atop 
a^{1\over 2}q^{-{1\over 2}},-a^{1\over 2}q^{-{1\over 2}},
aqb^{-1}};q,zqb^{-1}
\right).
\end{eqnarray}
\end{lem}\label{Lem1}
\proof
\begin{eqnarray*}
{\rm LHS}&=&1+\sum_{n=1}^\infty
{(a,b;q)_n\over (aqb^{-1},q;q)_n}(qb^{-1})^n z^n-
\sum_{n=0}^\infty
{(a,b;q)_n\over (aqb^{-1},q;q)_n}(qb^{-1})^nz^{n+1}\\
&=&1+\sum_{n=1}^\infty
{(a,b;q)_n\over (aqb^{-1},q;q)_n}
{(1-bq^{-1})(1-aq^{2n-1})\over
(1-aq^{n-1})(1-bq^{n-1})}
(zqb^{-1})^n ={\rm RHS}.
\end{eqnarray*}
\qed\\
\noindent
 We will see this identity will plays an important role in what follows.
 \medskip

\subsection{monomial basis representation of $I(\alpha)$}

We represent the action of the integral transformation $I(\alpha)$
in terms of the monomial basis $\{(\zeta_2/\zeta_1)^i|i=0,1,2,\cdots\}$ 
of ${\cal F}_2$. To this end, first we need 
the Laurent series expansion of the function
\begin{eqnarray}
&&g(\xi/\zeta)
{\Theta_{q}(\alpha  q^{1\over 2}t^{-{1\over 2}}\zeta/\xi)\over 
\Theta_{q}(q^{1\over 2}t^{-{1\over 2}}\zeta/\xi)}.
\end{eqnarray}
This is given by Ramanujan's summation formula for the ${}_1\psi_1$ series.
(Eq. (5.2.1) of Gasper and Rahman \cite{GR}, hereafter referred to as GR).
\begin{lem}\label{Ramanujan}
We have
\begin{eqnarray}
&&g(\xi/\zeta)
{\Theta_{q}(\alpha q^{1\over 2}t^{-{1\over 2}}\zeta/\xi)\over 
\Theta_{q}(q^{1\over 2}t^{-{1\over 2}}\zeta/\xi)}\label{gTheta}\\
&=&
{ 
(\alpha q t^{-1};q)_\infty
\over 
(q t^{-1};q)_\infty }
{ 
(\alpha^{-1} q;q)_\infty
\over 
(q;q)_\infty }
\sum_{m\in{\bf Z}}
{(\alpha;q)_m \over (\alpha q t^{-1};q)_m}
\left(q^{1\over 2} t^{-{1\over 2}} \zeta/\xi\right)^m ,\nonumber
\end{eqnarray}
for $|q^{1\over 2}t^{-{1\over 2}}|<|\zeta/\xi|<|q^{-{1\over 2}}t^{{1\over 2}}|$.
\end{lem}

Let us represent the action of $I(\alpha)$ 
on a formal power series of the form
\begin{eqnarray}
f(\zeta_1,\zeta_2)&=&
\sum_{j=0}^\infty f_j (\zeta_2/\zeta_1)^j, \label{c}
\end{eqnarray}
in terms of the basic hypergeometric series ${}_2\phi_1$.

\begin{prop}\label{Prop1}
Let $f(\zeta_1,\zeta_2)$ be as above.
Then we have
\begin{eqnarray}
&&I(\alpha)  f(\zeta_1,\zeta_2) \nonumber\\
&=&
(1-\zeta)
{(q t^{-1}\zeta;q)_\infty \over
 (t\zeta;q)_\infty }\\
&\times &
\sum_{k=0}^\infty 
{}_2\phi_1\left(
{ \alpha^{-1}s_1 q^{k }t,t\atop
 \alpha^{-1} s_1 q^{k+1}}; q,q t^{-1}\zeta
\right)
{}_2\phi_1\left(
{ \alpha s_2^{-1} q^{k },t\atop
 \alpha s_2^{-1} q^{k+1}t^{-1}}; q,q t^{-1}\zeta
\right)  \nonumber\\
&\times&
\lambda_k(\alpha) f_k \zeta^k, \nonumber
\end{eqnarray}
where $\zeta=\zeta_2/\zeta_1$.
\end{prop}

\proof 
By the $q$-binomial theorem (Eq. (1.3.2) of GR \cite{GR}), we have
\begin{eqnarray}
g(\zeta)=\sum_{j=0}^\infty g_j \zeta^j,
\qquad
g_j=
{(t;q)_i
\over 
(q;q)_i}
(q^{1\over 2}t^{-{1\over 2}})^i. \label{g_j}
\end{eqnarray}
Thus
\begin{eqnarray*}
&&{\rm LHS}\\
&=&
h(\zeta_2/\zeta_1)\oint{d\xi_1\over 2\pi i \xi_1}\oint{d\xi_2\over 2\pi i \xi_2}
\sum_{m,n\in{\bf Z}}
{(\alpha s_1^{-1};q)_m \over
(\alpha s_1^{-1}q t^{-1};q)_m }
(q^{1\over 2}t^{-{1\over 2}}\zeta_1/\xi_1)^m\\
&&\times
{(\alpha s_2^{-1};q)_n \over
(\alpha s_2^{-1}q t^{-1};q)_n }
(q^{1\over 2}t^{-{1\over 2}}\zeta_2/\xi_2)^n\\
&&\times
\sum_{i,j,k\geq 0} g_ig_j f_k (\xi_2/\zeta_1)^i
(\zeta_2/\xi_1)^j(\xi_2/\xi_1)^k\\
&=&h(\zeta_2/\zeta_1)
\sum_{i,j,k\geq 0} 
{(\alpha s_1^{-1};q)_{-j-k} \over
(\alpha s_1^{-1}q t^{-1};q)_{-j-k} } g_j
(q^{1\over 2}t^{-{1\over 2}})^{-j}(\zeta_2/\zeta_1)^j\\
&&\times
{(\alpha s_2^{-1};q)_{i+k} \over
(\alpha s_2^{-1}q t^{-1};q)_{i+k} }g_i
(q^{1\over 2}t^{-{1\over 2}} )^{i}(\zeta_2/\zeta_1)^i
\times f_k(\zeta_2/\zeta_1)^k \\
&=&h(\zeta_2/\zeta_1)
\sum_{k= 0}^\infty \sum_{j= 0}^\infty 
{(\alpha^{-1}s_1 q^{k}t,t;q)_{j} \over
(\alpha^{-1} s_1 q^{k+1},q;q)_{j} } 
(q t^{-1}\zeta_2/\zeta_1)^j\\
&&\times
\sum_{i= 0}^\infty 
{(\alpha s_2^{-1} q^{k},t;q)_{i} \over
(\alpha s_2^{-1} q^{k+1}t^{-1},q^{k};q)_{i} }
(q t^{-1} \zeta_2/\zeta_1)^i \\
&&\times
{(\alpha^{-1}s_1 t,\alpha s_2^{-1};q)_k \over
(\alpha^{-1} s_1 q,\alpha s_2^{-1}q t^{-1};q)_k} 
(q t^{-1} \zeta_2/\zeta_1)^k f_k\\
&=&{\rm RHS}.
\end{eqnarray*}
\qed

To obtain a more useful series expression for $I(\alpha)$ from Proposition \ref{Prop1},
we need a lemma.
\begin{lem}\label{Lem3}
We have
\begin{eqnarray}
&&{(z q/b;q)_\infty \over (bz;q)_\infty}
{~}_2\phi_1\left(
{ a,b \atop aq/b} ; q,zq/b
\right)
{~}_2\phi_1\left(
{ c,b \atop cq/b} ; q,zq/b
\right) \\
&=&
\sum_{k=0}^\infty 
{(cq/b^2;q)_k (q/b;q)_k \over 
(cq/b;q)_k (q;q)_k }b^k z^k
{~}_4\phi_3\left(
{ q^{-k},b,b q^{-k}/c,a \atop
 bq^{-k}, b^2q^{-k}/c,aq/b } ; q,q
\right) \nonumber.
\end{eqnarray}
\end{lem}

\proof
By using Heine's transformation (Eq. (1.4.3) of GR \cite{GR})
\begin{eqnarray*}
&&{\rm LHS}\\
&=&
{~}_2\phi_1\left(
{ a,b \atop aq/b} ; q,zq/b
\right)
{~}_2\phi_1\left(
{ cq/b^2 ,q/b  \atop cq/b} ; q,bz
\right) \\
&=&
\sum_{m,n=0}^\infty
{(a;q)_m (b;q)_m\over (aq/b;q)_m(q;q)_m}
{(cq/b^2;q)_n (q/b;q)_n\over (cq/b;q)_n(q;q)_n}
(zq/b)^m(bz)^n\\
&=&
\sum_{k=0}^\infty \sum_{m=0}^k
{(a;q)_m (b;q)_m\over (aq/b;q)_m(q;q)_m}
{(cq/b^2;q)_{k-m} (q/b;q)_{k-m}\over (cq/b;q)_{k-m}(q;q)_{k-m}}
(zq/b)^m(bz)^{k-m}\\
&=&
\sum_{k=0}^\infty 
{(cq/b^2;q)_{k} (q/b;q)_{k}\over (cq/b;q)_{k}(q;q)_{k}} (bz)^k\\
&&\times
\sum_{m=0}^k
{(a;q)_m (b;q)_m\over (aq/b;q)_m(q;q)_m}
{(bq^{-k}/c;q)_{k} (q^{-k};q)_{k}\over 
(b^2q^{-k}/c;q)_{k}(b q^{-k};q)_{k}}
q^m\\
&=&{\rm RHS}.
\end{eqnarray*}
\qed

Using Lemma \ref{Lem3},
we can express the action of the operator
$(1-\zeta_2/\zeta_1)^{-1}I(\alpha)$ 
by using a terminating and balanced ${}_4\phi_3$ series.
\begin{prop}\label{Prop3}
\begin{eqnarray}
{1\over 1-\zeta_2/\zeta_1 }I(\alpha)  f(\zeta_1,\zeta_2)
&=& \sum_{i=0}^\infty \sum_{j=0}^i (\zeta_2/ \zeta_1)^i 
e_{ij} f_j  , \label{eq1}
\end{eqnarray}
where
\begin{eqnarray}
e_{ij}&=&
{(q t^{-1};q)_{i-j}
(\alpha^{-1}s_1 q^{j+1}t^{-1};q)_{i-j}
\over 
(q ;q)_{i-j}
(\alpha^{-1} s_1 q^{j+1};q)_{i-j} }
t^{i-j} \lambda_j(\alpha) \\
&&\times
{}_4\phi_3\left(
{ q^{{-i+j}},t,\alpha s_1^{-1}q^{-{i}},\alpha s_2^{-1}q^{j}
\atop
 q^{{-i+j}}t,\alpha  s_1^{-1}q^{-{i}}t,\alpha  s_2^{-1}q^{j+1}t^{-1}}
 ; 
q,q\nonumber
\right).
\end{eqnarray}
\end{prop}

Note that, if we work with the monomial basis 
$\{(\zeta_2/\zeta_1)^i|i=0,1,2,\cdots\}$
and identifying $(\zeta_2/\zeta_1)^k$ with the column unit vector
$ {}^t(0,\cdots,0,1,0,\cdots)$ having $1$ at the $k$-th place,
we have an infinite lower triangular matrix representation of 
the integral operator $I(\alpha)$,
\begin{eqnarray}
{1\over 1-\zeta_2/\zeta_1 }I(\alpha)
=
\left(
\begin{array}{ccccc}
1& & & & \\
e_{10}& 1& & & \\
e_{20}& e_{21} & 1& & \\
e_{30}& e_{31}& e_{32}& 1& \\
\vdots & \vdots & & \ddots&\ddots
\end{array}
\right). 
\end{eqnarray}
Because of this lower triangular nature, 
the action of $I(\alpha)$ on ${\cal F}_n$ is well defined.

\proof
\begin{eqnarray*}
&&\mbox{LHS of Eq.(\ref{eq1})}\\
&=&
{(q t^{-1}\zeta;q)_\infty \over
 (t\zeta;q)_\infty}
\sum_{j=0}^\infty 
{}_2\phi_1\left(
{ \alpha^{-1}s_1 q^{j}t,t \atop \alpha^{-1}s_1 q^{j+1}} ; 
q,q t^{-1}\zeta
\right)\\
&&\times
{}_2\phi_1\left(
{ \alpha s_2^{-1} q^{j},t \atop \alpha s_2^{-1}q^{j+1}t^{-1}} ; 
q,q t^{-1}\zeta
\right) \lambda_j(\alpha)f_j \zeta^j\\
&=&
\sum_{j=0}^\infty 
\sum_{k=0}^\infty
{(q t^{-1};q)_k 
(\alpha^{-1}s_1 q^{j+1}t^{-1};q)_k 
\over 
(q;q)_k 
(\alpha^{-1} s_1 q^{j+1};q)_k }
(t\zeta)^k\lambda_j(\alpha)f_j \zeta^j\\
&&\times
{}_4\phi_3\left(
{ q^{-{k}},t,\alpha s_1^{-1} q^{-{j+k}},\alpha s_2^{-1}q^{j}
\atop
 q^{-{k}}t,\alpha  s_1^{-1} q^{-{j+k}}t,\alpha s_2^{-1} q^{j+1}t^{-1}}
 ; 
q,q
\right)\\
&=&
\sum_{i=0}^\infty 
\sum_{j=0}^i
{(q t^{-1};q)_{i-j}
(\alpha^{-1}s_1 q^{j+1}t^{-1};q)_{i-j}
\over 
(q ;q)_{i-j}
(\alpha^{-1} s_1 q^{j+1};q)_{i-j} }
t^{i-j} \lambda_j(\alpha)f_j \zeta^i\\
&&\times
{}_4\phi_3\left(
{ q^{{-i+j}},t,\alpha s_1^{-1}q^{-{i}},\alpha s_2^{-1}q^{j}
\atop
 q^{{-i+j}}t,\alpha  s_1^{-1}q^{-{i}}t,\alpha  s_2^{-1}q^{j+1}t^{-1}}
 ; 
q,q
\right)\\
&=&\mbox{RHS of Eq.(\ref{eq1})}.
\end{eqnarray*}
\qed

\subsection{properties of the coefficients of the functions $f_j(\zeta_1,\zeta_2)$}

Nextly, we need to study some properties of 
the coefficients of the functions $f_j(\zeta_1,\zeta_2)$ 
defined by Eq.(\ref{f_j})
in some detail.
the coefficients of $f_j(\zeta_1,\zeta_2)$ 
\begin{eqnarray}
f_j(\zeta_1,\zeta_2)=\sum_{i=0}^\infty c_{ij}(\zeta_2/\zeta_1)^i,
\end{eqnarray}
are given by
\begin{eqnarray}
c_{ij}&=&
\left\{
\begin{array}{ll}
\displaystyle
{
(s q^{2j}t^{-1}, s^{1\over 2}q^{j+1}t^{-{1\over 2}},
 -s^{1\over 2}q^{j+1}t^{-{1\over 2}},t^{-1};q)_{i-j}
\over 
(s q^{2j+1}, s^{1\over 2}q^{j}t^{-{1\over 2}},
 -s^{1\over 2}q^{j}t^{-{1\over 2}},q;q)_{i-j}}
t^{i-j} & (i\geq j) ,\\[4mm]
0 & (i<j) ,
\end{array}
\right.
\end{eqnarray}
where we have denoted $s=s_1/s_2$ for short.
If we identify the functions $f_j(\zeta_1,\zeta_2)$'s with the
column vectors ${}^t (c_{0j},c_{1j},c_{2j},\cdots)$
and use $\{f_j(\zeta_1,\zeta_2)|j=0,1,2\cdots\}$ as our basis of ${\cal F}_2$,
we may have another matrix representation of $I(\alpha)$.
With this basis, Theorem \ref{Thm2} is recast as 
\begin{eqnarray}
I(\alpha) C
=
C \Lambda(\alpha),
\end{eqnarray}
where
\begin{eqnarray}
C&=&\left(
\begin{array}{ccccc}
1& & & & \\
c_{10}& 1& & & \\
c_{20}& c_{21} & 1& & \\
c_{30}& c_{31}& c_{32}& 1& \\
\vdots & \vdots & & \ddots&\ddots
\end{array}
\right),\label{MatC}\\
\Lambda(\alpha)&=&{\rm diag}
(\lambda_0(\alpha),\lambda_1(\alpha),\lambda_2(\alpha),\cdots).
\end{eqnarray}

We find that
all the entries of $C^{-1}$ can be factorized.
\begin{prop}
The inverse of C is written as
\begin{eqnarray}
C^{-1}&=&\left(
\begin{array}{ccccc}
1& & & & \\
d_{10}& 1& & & \\
d_{20}& d_{21} & 1& & \\
d_{30}& d_{31}& d_{32}& 1& \\
\vdots & \vdots & & \ddots&\ddots
\end{array}
\right),
\end{eqnarray}
where
\begin{eqnarray}
d_{ij}&=&
\left\{
\begin{array}{ll}
\displaystyle
{
(s q^{i+j+1}t^{-1},t;q)_{i-j}
\over 
(s q^{i+j},q;q)_{i-j}} & (i\geq j) ,\\[4mm]
0 & (i<j) .
\end{array}
\right.
\end{eqnarray}

\end{prop}

\proof
It suffices to show
\begin{eqnarray*}
\sum_{k=j}^i d_{ik}c_{kj} =\delta_{i,j},
\end{eqnarray*}
for $i\geq j$.
Rewrite the left hand side as
$d_{i,j}\sum_{k=0}^{i-j}{d_{i,k+j} \over d_{i,j}}c_{k+j,j}$, and
we realize this summation 
is a terminating very-well-poised ${}_6\phi_5$ series, since 
\begin{eqnarray*}
{d_{i,k+j} \over d_{i,j}}
&=&{
(q^{-i+j},s q^{i+j};q)_k
\over 
(q^{-i+j+1}t^{-1},s q^{i+j+1}t^{-1};q)_k}
(q t^{-1})^k,\\
c_{k+j,j}&=&
{
(s q^{2j} t^{-1}, s^{1\over 2} q^{j+1}t^{-{1\over 2}},
 - s^{1\over 2}q^{j+1}t^{-{1\over 2}},t^{-1};q)_k
\over 
(s q^{2j+1},  s^{1\over 2}q^{j}t^{-{1\over 2}},
 - s^{1\over 2}q^{j}t^{-{1\over 2}},q;q)_k}
t^k.
\end{eqnarray*}
One then finds that 
this ${}_6\phi_5$ series can be summed by 
using Eq. (2.4.2) of GR \cite{GR}. Thus we have
\begin{eqnarray*}
&&d_{i,j}\sum_{k=0}^{i-j}{d_{i,k+j} \over d_{i,j}}c_{k+j,j}\\
&=&d_{i,j}
{~}_6\phi_5
\left(
{
sq^{2j} t^{-1}, s^{1\over 2}q^{j+1}t^{-{1\over 2}},
-s^{1\over 2}q^{j+1}t^{-{1\over 2}}, t^{-1},s q^{i+j},
q^{-i+j}
\atop
s^{1\over 2}q^{j}t^{-{1\over 2}},
-s^{1\over 2}q^{j}t^{-{1\over 2}}, s q^{2j+1},
q^{-i+j+1}t^{-1},
s q^{i+j+1}t^{-1}};q,q
\right)\\
&=&d_{i,j} 
{
(s q^{2j+{1}}t^{-1},q^{-i+j+1};q)_{i-j}
\over
( s q^{2j+{1}},q^{-i+j+1}t^{-1};q)_{i-j}}\\
&=&\delta_{i,j}.
\end{eqnarray*}
\qed

Our next task is to study the matrix 
$C\Lambda(\alpha)C^{-1}$ and compare it with the coefficients given in
Proposition \ref{Prop3}. 
Since we havey realized $(1-\zeta_2/\zeta_1)^{-1}I(\alpha)$
instead of $I(\alpha)$ itself, 
we need to modify the function $f_j(\zeta_1,\zeta_2)$ accordingly.
Let us set
\begin{eqnarray}
f_j(\zeta_1,\zeta_2)&=&(1-\zeta_2/\zeta_1)\widetilde{f}_j(\zeta_1,\zeta_2),
\end{eqnarray}
where
\begin{eqnarray}
\widetilde{f}_j(\zeta_1,\zeta_2)&=&
(\zeta_2/\zeta_1)^j\times
{}_2\phi_1\left(
{  q^{2j+{1}} t^{-1}, q t^{-1}
\atop
q^{2j+{1}}}
;q, t\zeta_2/\zeta_1
\right).
\end{eqnarray}
Then Theorem \ref{Thm2} can be written as
\begin{eqnarray}
{1\over 1-\zeta_2/\zeta_1}I(\alpha) f_j(\zeta_1,\zeta_2) =\lambda_j (\alpha)
\widetilde{f}_j(\zeta_1,\zeta_2)
\qquad (j=0,1,2,\cdots). \label{If=lamf}
\end{eqnarray}
Setting
\begin{eqnarray}
\tilde{f}_j(\zeta_1,\zeta_2)=\sum_{i=0}^\infty \tilde{c}_{ij}(\zeta_2/\zeta_1)^i,
\end{eqnarray}
we have a matrix equation
\begin{eqnarray}
{1\over 1-\zeta_2/\zeta_1}I(\alpha)&=&
\widetilde{C}\Lambda(\alpha)C^{-1}, \label{EqThm}
\end{eqnarray}
where
\begin{eqnarray}
\widetilde{C}&=&(\widetilde{c}_{ij})_{0\leq i,j\leq \infty},\\
\widetilde{c}_{ij}&=&
\left\{
\begin{array}{ll}
\displaystyle
{
(s q^{2j+{1}}t^{-1},q t^{-1};q)_{i-j}
\over 
(s q^{2j+{1}},q;q)_{i-j}}
t^{i-j} & (i\geq j), \\[4mm]
0 & (i<j) ,
\end{array}
\right.
\end{eqnarray}
and $s=s_1/s_2$.
\medskip

Let us examine the matrix $\widetilde{C}C^{-1}$ before we
start to calculate $\widetilde{C}\Lambda(\alpha)C^{-1}$.
It is not necessary, but 
this way of calculating the matrix multiplication 
may reduce our task a little.
\begin{prop}\label{PropCtC}
We have
\begin{eqnarray}
\widetilde{C}C^{-1}=
\left(
\begin{array}{ccccc}
1& & & & \\
1& 1& & & \\
1& 1 & 1& & \\
1& 1& 1& 1& \\
\vdots&\vdots & & \ddots&\ddots
\end{array}
\right).
\end{eqnarray}
\end{prop}

\proof
By using
\begin{eqnarray*}
&&{(t;q)_{i-j-l} \over (q;q)_{i-j-l}}
=
{(t;q)_{i-j} \over (q;q)_{i-j}}
{(q^{-i+j};q)_{l} \over 
(q^{-i+j+1}t^{-1};q)_{l}}(q t^{-1})^l,\\[4mm]
&&
(s q^{2i-2l+{1}} t^{-1};q)_{l}
(s q^{i+j-l+1}t^{-1};q)_{i-j-l}
=
(s q^{i+j-l+1}t^{-1};q)_{i-j}\nonumber\\
&&\qquad =
(s q^{i+j+1}t^{-1};q)_{i-j}
{
( s^{-1}q^{-i-j}t;q)_l \over
( s^{-1}q^{-2i}t;q)_l}
q^{-l(i-j)},\\[4mm]
&&
(s q^{2i-2l+{1}};q)_{l}
(s q^{i+j-l};q)_{i-j-l}
=
{1-s q^{2i-{l}}\over 1-s q^{2i-2l}}
(s q^{i+j-l};q)_{i-j}\nonumber\\
&&\qquad =
{1-s q^{2i-{l}}\over 1-s q^{2i-2l}}
(s q^{i+j};q)_{i-j}
{
(s^{-1}q^{-i-j+1};q)_l\over
(s^{-1}q^{-2i+{1}};q)_l}q^{-l(i-j)},
\end{eqnarray*}
and
\begin{eqnarray*}
&&{1-s q^{2i-2l}\over 1-s q^{2i-{l}}}
=
{1-s^{-1} q^{-2i+2l}\over 1-s^{-1}q^{-2i}}
{1-s^{-1}q^{-2i}\over 1-s^{-1}q^{-2i+{l}}}q^{-l}\\
&&\qquad =
{(s^{-{1\over 2}}q^{-i+1},-s^{-{1\over 2}}q^{-i+1},
s^{-1}q^{-2i};q)_l
\over 
(s^{-{1\over 2}}q^{-i},-s^{-{1\over 2}}q^{-i},
s^{-1}q^{-2i+{1}};q)_l}q^{-l},\nonumber
\end{eqnarray*}
we have
\begin{eqnarray*}
&&\sum_{k=j}^i \widetilde{c}_{ik}d_{kj}
=\sum_{l=0}^{i-j}
\widetilde{c}_{i,i-l}d_{i-l,j}\\ 
&=&\sum_{l=0}^{i-j}
{(s q^{2i-2l+{1}}t^{-1},q t^{-1};q)_l
\over 
( s q^{2i-2l+{1}}, q;q)_l} t^l
{
(s q^{i+j-l+1}t^{-1},t;q)_{i-j-l}\over
(s q^{i+j-l},q;q)_{i-j-l} } \\
&=&
d_{ij}\;{}_6\phi_5
\left(
{
s^{-1}q^{-2i},s^{-{1\over 2}}q^{-i+1},-s^{-{1\over 2}}q^{-i+1},
q t^{-1},
 s^{-1}q^{-i-j}t,q^{-i+j}
\atop
s^{-{1\over 2}}q^{-i},-s^{-{1\over 2}}q^{-i},
 s^{-1}q^{-2i}t,q^{-i+j+1}t^{-1},s^{-1}q^{-i-j+1}}
;q,1
\right)\\
&=&
d_{ij}
{
(s^{-1} q^{-2i+{1}},q^{-i+j};q)_{i-j}
\over
(s^{-1} q^{-2i }t ,q^{-i+j+1}t^{-1};q)_{i-j}
}\\
&=&1.
\end{eqnarray*}
Here, we have used Eq. (2.4.2) of GR \cite{GR}.
\qed

Now we proceed to studying the matrix elements of
 $\widetilde{C}\Lambda(\alpha)C^{-1}$. 
\begin{prop}
The nonzero entries for $\widetilde{C}\Lambda(\alpha)C^{-1}$ are 
are given by
\begin{eqnarray}
&&(\widetilde{C}\Lambda(\alpha)C^{-1})_{ij} \nonumber\\
&=&
\lambda_i(\alpha)d_{ij} \label{CtLamC}\\
&\times&
{}_8W_7\left( 
s^{-1}q^{-2i};
q t^{-1},s^{-1} q^{-i-j}t,
\alpha s_1^{-1}q^{-i},\alpha^{-1} s_2 q^{-i}t,q^{-i+j}
; q,q t^{-1}
\right), \nonumber
\end{eqnarray}
for $i\geq j$, and $(\widetilde{C}\Lambda(\alpha)C^{-1})_{ij}=0$
if $i<j$.
\end{prop}
Note that here we have used the compact notation
\begin{eqnarray}
&&{}_8W_7( a_1;a_4,a_5,\cdots,a_{r+1};q,z)\\
&=&
{}_{r+1}\phi_r\left(
{a_1,q a_1^{1\over 2},-q a_1^{1\over 2},a_4,\cdots,a_{r+1}
\atop
a_1^{1\over 2},-a_1^{1\over 2},qa_1/a_4,\cdots,q a_1/a_{r+1}
};q,z
\right),\nonumber
\end{eqnarray}
for the very-well-poised ${}_{r+1}\phi_r$ series.

\proof
From the identity
\begin{eqnarray*}
&&\lambda_{i-l}(\alpha)\\
&=&
{(\alpha s_2^{-1},\alpha^{-1} s_1 t;q)_{i-l}
\over 
(\alpha s_2^{-1}q t^{-1},\alpha^{-1}s_1q;q)_{i-l}}
(q t^{-1})^{i-l}\\
&=&
{(\alpha s_2^{-1},\alpha^{-1}s_1t;q)_{i}
\over 
(\alpha s_2^{-1}q t^{-1},\alpha^{-1}s_1q;q)_{i}}
(q t^{-1})^{i}
{(\alpha^{-1} s_2 q^{-i}t,\alpha s_1^{-1}q^{-i};q)_{l}
\over 
(\alpha^{-1} s_2 q^{-i+1},\alpha s_1^{-1}q^{-i+1}t^{-1};q)_{l}}
(q t^{-1})^{l}\\
&=&
\lambda_i(\alpha)
{(\alpha^{-1} s_2 q^{-i}t,\alpha s_1^{-1}q^{-i};q)_{l}
\over 
(\alpha^{-1} s_2 q^{-i+1},\alpha s_1^{-1}q^{-i+1}t^{-1};q)_{l}}
(q t^{-1})^{l},
\end{eqnarray*}
and the calculation given in the proof of Proposition \ref{PropCtC}, we have
\begin{eqnarray*}
\mbox{LHS of Eq.(\ref{CtLamC})}
&=&\sum_{k=j}^i \widetilde{c}_{ik}\lambda_k(\alpha)d_{kj}
= \sum_{l=0}^{i-j}
\widetilde{c}_{i,i-l}
 \lambda_{i-l}(\alpha) d_{i-l,j}\\ 
&=&\mbox{RHS of Eq.(\ref{CtLamC})}.
\end{eqnarray*}

\qed

Then, one can apply Watson's formula (Eq. (2.5.1) of GR \cite{GR})
to transform the above terminating ${}_8W_7$
series
into a terminating balanced ${}_4\phi_3$ series.
\begin{prop}
We have
\begin{eqnarray}
&&(\widetilde{C}\Lambda(\alpha)C^{-1})_{ij} 
=
\lambda_i(\alpha)d_{ij} 
{
(s^{-1}q^{-2i+{1}}, q t^{-1};q)_{i-j} \over
(\alpha^{-1} s_2 q^{{-i+1}}, \alpha s_1^{-1}q^{-i+1}t^{-1};q)_{i-j}
}\\
&&\qquad \times
{}_4\phi_3\left(
{ q^{{-i+j}},\alpha s_1^{-1}q^{-{i}},\alpha^{-1} s_2 q^{-{i}}t,
q^{{-i+j}}
\atop
 q^{{-i+j}}t, q^{{-i+j+1}}t^{-1},
 s^{-1}q^{-2i}t}
 ; 
q,q
\right),\nonumber
\end{eqnarray}
for $i\geq j$.
\end{prop}

Finally, by using  Sears' transformation 
(Eq. (2.10.4) of GR \cite{GR}),
we see that the coefficients  $(\widetilde{C}\Lambda(\alpha)C^{-1})_{ij}$ 
and $e_{ij}$ are the same.
\begin{prop}\label{Prop4}
We have
\begin{eqnarray}
&&(\widetilde{C}\Lambda(\alpha)C^{-1})_{ij} 
=e_{ij}
\end{eqnarray}
for $i\geq j$.
\end{prop}

\proof
By using
\begin{eqnarray*}
&&{
(t;q)_{i-j}\over
(q^{-i+j+1}t^{-1};q)_{i-j}
}
=(-1)^{i-j}t^{i-j}q^{(i-j)(i-j-1)\over 2},\\
&&{
(s^{-1}q^{-2i+{1}};q)_{i-j}\over
(s q^{i+j};q)_{i-j}
}
=(-1)^{i-j}s^{-i+j}q^{-{(i+j)(i-j)}}q^{-{(i-j)(i-j-1)\over 2}},\\
&&{
(s q^{{i+j+1}}t^{-1};q)_{i-j}\over
(s^{-1}q^{-2i}t;q)_{i-j}
}
=(-1)^{i-j}s^{i-j}t^{-i+j}q^{(i+j+1)(i-j)}q^{{(i-j)(i-j-1)\over 2}},
\end{eqnarray*}

\begin{eqnarray*}
&&{
(\alpha^{-1}s_2 q^{{-i}}t;q)_{i-j}\over
(\alpha^{-1}s_2 q^{-i+1};q)_{i-j}
}
=
{
(\alpha s_2^{-1}q t^{-1};q)_{i}\over
(\alpha s_2^{-1};q)_{i}
}
{
(\alpha s_2^{-1};q)_{j}\over
(\alpha s_2^{-1}q t^{-1};q)_{j}
}(q t^{-1})^{-i+j},\\
&&{
(\alpha^{-1}s_1 q^{{j+1}};q)_{i-j}\over
(\alpha s_1^{-1} q^{-i+1}t^{-1};q)_{i-j}
}
=
{
(\alpha^{-1} s_1q ;q)_{i}\over
(\alpha^{-1} s_1t;q)_{i}
}
{
(\alpha^{-1} s_1t;q)_{j}\over
(\alpha^{-1} s_1q;q)_{j}
}\\
&&\qquad\qquad\times
(-1)^{-i+j}\alpha^{-i+j}s_1^{i-j}t^{i-j}q^{{i(i-1)\over 2}-{j(j-1)\over 2}},
\end{eqnarray*}
we have
\begin{eqnarray*}
&&(\widetilde{C}\Lambda(\alpha)C^{-1})_{ij} \\
&=&
\lambda_i(\alpha)d_{ij} 
{
(s^{-1}q^{-2i+{1}}, q t^{-1};q)_{i-j} \over
(\alpha^{-1} s_2 q^{{-i+1}}, \alpha s_1^{-1}q^{-i+1}t^{-1};q)_{i-j}
}\\
&&\qquad\times
{
(\alpha^{-1} s_2 q^{{-i}}t, \alpha^{-1} s_1 q^{j}t^{-1};q)_{i-j}
\over
(s^{-1} q^{-2i} t,q^{-i+j+1}t^{-1};q)_{i-j}
}(\alpha s_1^{-1} q^{-i})^{i-j}\\
&&\qquad \times
{}_4\phi_3\left(
{ q^{{-i+j}},\alpha s_1^{-1}q^{-{i}},\alpha s_2^{-1}q^{{j}},t
\atop
 q^{{-i+j}}t, \alpha s_1^{-1} q^{-{i}}t, \alpha s_2^{-1} q^{{j+1}}t^{-1}}, 
 ; 
q,q
\right)\\
 &=&
 {(q t^{-1};q)_{i-j}
(\alpha^{-1}s_1 q^{j+1}t^{-1};q)_{i-j}
\over 
(q ;q)_{i-j}
(\alpha^{-1} s_1 q^{j+1};q)_{i-j} }
t^{i-j} \lambda_j(\alpha) \\
&&\times
{}_4\phi_3\left(
{ q^{{-i+j}},t,\alpha s_1^{-1}q^{-{i}},\alpha s_2^{-1}q^{j}
\atop
 q^{{-i+j}}t,\alpha  s_1^{-1}q^{-{i}}t,\alpha  s_2^{-1}q^{j+1}t^{-1}}
 ; 
q,q \nonumber
\right)=e_{ij}.
\end{eqnarray*}
\qed
\medskip

We obtain Eq.(\ref{If=lamf}) from Proposition \ref{Prop3} and
Proposition \ref{Prop4}. Hence
we have completed the proof of Theorem \ref{Thm2}.

\section{Properties of the Eigenfunctions}
We examine some basic properties of the eigenfunctions of $I(\alpha)$ for general $n$.
An explicit formua for all the eigenfunctions of $I(\alpha)$ for $n=3$
is conjectured, and a partial result is given for the case of $n=4$. 
In these explicit formulas,
no dependence on the spectral parameter $\alpha$ is observed.
This supports our Conjecture \ref{Con1}.

\subsection{existence of the eigenfunctions}
Let us study the existence of the eigenfunctions
which form a basis of ${\cal F}_n$. 
\begin{prop}\label{existing}
Let the parameters $(s_1,s_2,\cdots,s_n)$, $\alpha$, $q$ and $t$ be generic.
Let $j_1,j_2,\cdots,j_{n-1}$ be nonnegative integers.
In the space ${\cal F}_n$, there exist a unique solution to the equation
\begin{eqnarray}
&&I(\alpha) f_{j_1,j_2,\cdots,j_{n-1}}(\zeta_1,\cdots,\zeta_n)\\
&=&
 \lambda_{j_1,j_2,\cdots,j_{n-1}} (\alpha) 
 f_{j_1,j_2,\cdots,j_{n-1}}(\zeta_1,\cdots,\zeta_n),
 \nonumber
\end{eqnarray}
with the conditions
\begin{eqnarray}
&&
f_{j_1,j_2,\cdots,j_{n-1}}(\zeta_1,\cdots,\zeta_n)\\
&=&
\sum_{i_1\geq j_1, \cdots,i_{n-1}\geq j_{n-1}}^\infty
c_{i_1,i_2,\cdots,i_{n-1}}\left(\zeta_2\over \zeta_1\right)^{i_1}
\left(\zeta_3\over \zeta_2\right)^{i_2}
\cdots \left(\zeta_n\over \zeta_{n-1}\right)^{i_{n-1}},\nonumber
\end{eqnarray}
and $c_{j_1,j_2,\cdots,j_{n-1}}=1$,
if and only if 
\begin{eqnarray}
\lambda_{j_1,j_2,\cdots,j_{n-1}}(\alpha)
&=&
\prod_{i=1}^{n}
{(\alpha s_i^{-1};q)_{j_{i-1}-j_i} \over 
  (\alpha s_i^{-1} q t^{-1};q)_{j_{i-1}-j_i}},\label{eigenvalues}
\end{eqnarray}
is satisfied.
Here $j_0=0$ and $j_{n}=0$ are assumed.
\end{prop}

\proof
Working with the dominance order, or the lexicographic order in the monomial basis
\begin{eqnarray}
1<{\zeta_2\over \zeta_1}<\left({\zeta_2\over \zeta_1}\right)^2<\cdots<
{\zeta_3\over \zeta_2}<{\zeta_2\over \zeta_1}{\zeta_3\over \zeta_2}<
\left({\zeta_2\over \zeta_1}\right)^2{\zeta_3\over \zeta_2}<\cdots,\nonumber
\end{eqnarray}
we see the representation of the integral operation $I(\alpha)$ becomes
infinite lower triangular matrix. The explicit form of the action of $I(\alpha)$ 
on a monomial is 
obtained by Lemma \ref{Ramanujan}.
It reads
\begin{eqnarray*}
&&I(\alpha) \left(\zeta_2\over \zeta_1\right)^{j_1}
\left(\zeta_3\over \zeta_2\right)^{j_2}
\cdots \left(\zeta_n\over \zeta_{n-1}\right)^{j_{n-1}}\\
&=&
\left(\zeta_2\over \zeta_1\right)^{j_1}
\left(\zeta_3\over \zeta_2\right)^{j_2}
\cdots \left(\zeta_n\over \zeta_{n-1}\right)^{j_{n-1}}\prod_{i<j}h(\zeta_j/\zeta_i)\\
&\times&
\sum_{k_{l,m}\geq 0\atop ( 1 \leq l\neq m\leq n)}
\prod_{r=1}^n
\mu(\alpha s_r^{-1}; k_{1,r}+\cdots +k_{r-1,r}-k_{r+1,r}-\cdots-k_{n,r}+j_{r-1}-j_r)\\
&\times&
\prod_{l<m} g_{k_{l,m}}
\left(\zeta_m\over \zeta_l\right)^{k_{l,m}}
\prod_{l>m} g_{k_{l,m}}
\left(\zeta_l\over \zeta_m\right)^{k_{l,m}},
\end{eqnarray*}
where we have used the notations Eq.(\ref{g_j}) and
\begin{eqnarray*}
\mu(\alpha;k)={(\alpha;q)_k \over 
(\alpha q t^{-1};q)_k}
(q^{1\over 2}t^{-{1\over 2}})^k.
\end{eqnarray*}
Thus the lower triangularity is explicitly seen.
It is easy to see that 
the diagonal elements are given by Eq.(\ref{eigenvalues}),
by setting $k_{l,m}=0$ and forgetting the factor $\prod_{i<j} h(\zeta_j/\zeta_i)$ 
in the above expression. 
Since all the diagonal entries are distinct if the parameters are generic, there is
no obstruction for the construction of the eigenfunction.
\qed

By noting that 
\begin{eqnarray*}
\mu(\alpha;k+l)=\mu(\alpha;k)\mu(\alpha q^{k},l),
\end{eqnarray*}
it can be easily seen that all the eigenfunctions are related by 
shifting the parameters $s_i$ suitably.

\begin{prop}\label{eigen-rel}
The eigenfunctions of $I(\alpha)$ 
satisfy
\begin{eqnarray}
&&f_{j_1,j_2,\cdots,j_{n-1}}(\zeta_1,\cdots,\zeta_n)\\
&=&
\prod_{i=1}^{n}
 \zeta_i^{j_{i-1}-j_i}
(T_{q,s_i})^{-j_{i-1}+j_i}\cdot
f_{0,0,\cdots,0}(\zeta_1,\cdots,\zeta_n).\nonumber
\end{eqnarray}
Here, $j_0=0,j_n=0$ are assumed, and the shift operator 
$T_{q,s_i}$ is defined by
\begin{eqnarray}
T_{q,s_i} \cdot g(s_1,\cdots,s_n)=
g(s_1,\cdots,q s_i,\cdots s_n).
\end{eqnarray}

\end{prop}

We have constructed 
the set of eigenfunctions 
$\{f_{j_1,j_2,\cdots,j_{n-1}}(\zeta_1,\cdots,\zeta_n)\}$
which forms a basis of ${\cal F}_n$. 
Note that, we have a good structure
in the eigenfunctions
associated to the $A_{n-1}$ root system. Namely, 
the relation between the leading term and 
the shifts in the parameters $s_i$
are nicely organized 
by the $A_{n-1}$ roots
$\alpha_1=(-1,1,0,0,\cdots)$, 
$\alpha_2=(0,-1,1,0,\cdots)$ and so on.
\medskip

\subsection{conjecture for $n=3$}

Let us examine the eigenfunctions for $n=3$ and argue that
our  Conjecture \ref{Con1} seems true for $n=3$. 
We have not obtained a good enough understanding of $I(\alpha)$
to be able to derive the eigenfunctions for $n\geq 3$ in a rigorous manner.
Fortunately, however, 
a conjectural form of the eigenfunctions can be obtained for $n=3$
by a brute force calculation up to certain degrees in $\zeta$-variables.

\begin{con}\label{eigen_for_n=3}
The first eigenfunction of the integral transformation $I(\alpha)$ for $n=3$
is given by
\begin{eqnarray}
&&f_{0,0}(\zeta_1,\zeta_2,\zeta_3)\nonumber\\
&=&
\sum_{k=0}^\infty
{
(q t^{-1},q t^{-1},t,t;q)_k \over 
(q,q s_1/s_2,q s_2/s_3,
q s_1/s_3;q)_k}
 (q s_1/s_3)^k (\zeta_3/\zeta_1)^k\label{g-fun}\\
&&\times
\prod_{1\leq i<j\leq 3}
(1-\zeta_j/\zeta_i)
\cdot {}_2\phi_1\left(
{  q^{k+1} t^{-1}, q t^{-1}s_i/s_j ,
\atop
q^{{k+1}}s_i/s_j  };q, t\zeta_j/\zeta_i
\right).\nonumber
\end{eqnarray}
\end{con}
Note that all the other eigenfunctions can be obtained by Proposition \ref{eigen-rel}.

Since all the eigenfunctions given above 
clearly do not depend on the parameter $\alpha$, this strongly supports Conjecture \ref{Con1} for $n=3$.

\subsection{partial result for $n=4$}
For $n=4$, the study of the eigenfunctions becomes much more difficult, 
and we have not completely understood the structure of them yet.
So we give first several terms of the eigenfunctions to state our observation that
Conjecture \ref{Con1} seems true for $n=4$.
 
\begin{prop}\label{n=4-eigen}
The first eigenfunction 
$f_{0,0,0}(\zeta_1,\zeta_2,\zeta_3,\zeta_4)$
up to the powers
\begin{eqnarray}
\left({\zeta_2\over \zeta_1}\right)^{i_1}
\left({\zeta_3\over \zeta_2}\right)^{i_2}
\left ({\zeta_4\over \zeta_3}\right)^{i_3}\qquad
(0\leq  i_1\leq 2,0\leq  i_2\leq 2,0\leq  i_3\leq 2),
\end{eqnarray}
 is given by the following expression.
\begin{eqnarray}
&&f_{0,0,0}(\zeta_1,\zeta_2,\zeta_3,\zeta_4)\label{g-fun-4}\\
&=&
y_{0,0,0}+
y_{1,1,0}+y_{0,1,1}+y_{1,1,1}\nonumber \\
&& + \,  y_{2,2,0}+y_{1,2,1}+y_{0,2,2}
+
y_{2,2,1}+y_{1,2,2}+y_{2,2,2}+\cdots,\nonumber
\end{eqnarray}
where
\begin{eqnarray}
&&y_{0,0,0}= \varphi(0,0,0,0,0,0),\\
&&y_{1,1,0}
={\zeta_3\over \zeta_1} \left(q {s_1\over s_3} \right) 
{(q t^{-1})_1(q t^{-1})_1(t)_1(t)_1 \over
(q)_1(q s_{12})_1(q s_{23})_1(q s_{13})_1}
\varphi(1,1,0,1,0,0),\\
&&y_{0,1,1} 
= {\zeta_4\over \zeta_2} \left(q {s_2\over s_4} \right) 
{(q t^{-1})_1(q t^{-1})_1(t)_1(t)_1 \over
(q)_1(q s_{23})_1(q s_{34})_1(q s_{24})_1}
\varphi(0,1,1,0,1,0),\\
&&y_{1,1,1}\\
&& ={\zeta_4\over \zeta_1} \left(q {s_1\over s_4} \right) 
{(q t^{-1})_1(q t^{-1})_1(t)_1(t)_1 \over
(q)_1(q s_{12})_1(q s_{24})_1(q s_{14})_1}
\varphi(1,0,0,0,1,1)\nonumber \\
&&+ {\zeta_4\over \zeta_1} \left(q {s_1\over s_4} \right) 
{(q t^{-1})_1(q t^{-1})_1(t)_1(t)_1 \over
(q)_1(q s_{13})_1(q s_{34})_1(q s_{14})_1}
\varphi(0,0,1,1,0,1)\nonumber\\
&&+ {\zeta_4\over \zeta_1} \left(-q{s_1\over s_4} \right) 
{(q t^{-1})_1(q t^{-1})_1(q t^{-1})_1
(t)_1(t)_1(t)_1  (q^2 \,s_1s_2/s_3s_4)_1\over
(q)_1(q s_{12})_1(q s_{23})_1
(q s_{34})_1(q s_{13})_1(q s_{24})_1
(q s_{14})_1} \nonumber\\
&&\times
\varphi(1,1,1,1,1,1),\nonumber \\
&&y_{2,2,0}
= {\zeta_3\zeta_3\over \zeta_1\zeta_1} \left(q^2 {s_1s_1\over s_3s_3} \right) 
{(q t^{-1})_2(q t^{-1})_2(t)_2(t)_2 \over
(q)_2(q s_{12})_2(q s_{23})_2(q s_{13})_2}
\varphi(2,2,0,2,0,0),\\
&&y_{0,2,2} 
=  {\zeta_4\zeta_4\over \zeta_2\zeta_2} \left(q^2{s_2 s_2\over s_4 s_4} \right) 
{(q t^{-1})_2(q t^{-1})_2(t)_2(t)_2 \over
(q)_2(q s_{23})_2(q s_{34})_2(q s_{24})_2}
\varphi(0,2,2,0,2,0),\\
&&y_{1,2,1} \\
&&= {\zeta_3\zeta_4\over \zeta_1\zeta_2} \left(-q{s_1s_2\over s_3s_4} \right) 
{(q t^{-1})_1(q t^{-1})_1(q t^{-1})_1
(t)_1(t)_1(t)_1  (q^2 \,s_1s_3/s_2s_4)_1\over
(q)_1(q s_{12})_1(q s_{23})_1
(q s_{34})_1(q s_{13})_1(q s_{24})_1
(q s_{14})_1} \nonumber\\
&&\times
\varphi(1,1,1,1,1,1) \nonumber\\
&&+ {\zeta_3\zeta_4\over \zeta_1\zeta_2} \left(q{s_1s_2\over s_3s_4} \right) 
{(q t^{-1})_2(q t^{-1})_1(q t^{-1})_1
(t)_2(t)_1(t)_1  (q^2 \,s_{14})_1\over
(q)_1(q)_1(q s_{12})_1(q s_{23})_2
(q s_{34})_1(q s_{13})_1(q s_{24})_1
(q s_{14})_1}\nonumber\\
&&\times
\varphi(1,2,1,1,1,1),\nonumber\\
&&y_{2,2,1} \\
&&= {\zeta_3\zeta_4\over \zeta_1\zeta_1} \left(-q{s_1s_1\over s_3s_4} \right) 
{(q t^{-1})_1(q t^{-1})_1(q t^{-1})_1
(t)_1(t)_1(t)_1  (q^2 \,s_2s_3/s_1s_4)_1\over
(q)_1(q s_{12})_1(q s_{23})_1
(q s_{34})_1(q s_{13})_1(q s_{24})_1
(q s_{14})_1} \nonumber\\
&&\times
\varphi(1,1,1,1,1,1) \nonumber\\
&&+ {\zeta_3\zeta_4\over \zeta_1\zeta_1} \left(q{s_1s_1\over s_3s_4} \right) 
{(q t^{-1})_2(q t^{-1})_1(q t^{-1})_1
(t)_2(t)_1(t)_1  (q^2 \,s_{34})_1\over
(q)_1(q)_1(q s_{12})_2(q s_{23})_1
(q s_{34})_1(q s_{13})_1(q s_{24})_1
(q s_{14})_1}\nonumber\\
&&\times
\varphi(2,1,1,1,1,1)\nonumber\\
&&+ {\zeta_3\zeta_4\over \zeta_1\zeta_1} \left(q{s_1s_1\over s_3s_4} \right) 
{(q t^{-1})_2(q t^{-1})_1(q t^{-1})_1
(t)_2(t)_1(t)_1  (q^2 \,s_{24})_1\over
(q)_1(q)_1(q s_{12})_1(q s_{23})_1
(q s_{34})_1(q s_{13})_2(q s_{24})_1
(q s_{14})_1}\nonumber\\
&&\times
\varphi(1,1,1,2,1,1)\nonumber\\
&&+ {\zeta_3\zeta_4\over \zeta_1\zeta_1} \left(-q {s_1s_1\over s_3s_4} \right) 
{(q t^{-1})_2(q t^{-1})_2(q t^{-1})_1
(t)_2(t)_2(t)_1  (q^{3} s_1s_2/s_3s_4)_1\over
(q)_1(q)_1(q s_{12})_2(q s_{23})_2
(q s_{34})_1(q s_{13})_2(q s_{24})_1
(q s_{14})_1}\nonumber\\
&&\times
\varphi(2,2,1,2,1,1),\nonumber\\
&&y_{1,2,2} \\
&&= {\zeta_4\zeta_4\over \zeta_1\zeta_2} \left(-q{s_1s_2\over s_4s_4} \right) 
{(q t^{-1})_1(q t^{-1})_1(q t^{-1})_1
(t)_1(t)_1(t)_1  (q^2 \,s_1s_4/s_2s_3)_1\over
(q)_1(q s_{12})_1(q s_{23})_1
(q s_{34})_1(q s_{13})_1(q s_{24})_1
(q s_{14})_1} \nonumber\\
&&\times
\varphi(1,1,1,1,1,1) \nonumber\\
&&+ {\zeta_4\zeta_4\over \zeta_1\zeta_2} \left(q{s_1s_2\over s_4s_4} \right) 
{(q t^{-1})_2(q t^{-1})_1(q t^{-1})_1
(t)_2(t)_1(t)_1  (q^2 \,s_{12})_1\over
(q)_1(q)_1(q s_{12})_1(q s_{23})_1
(q s_{34})_2(q s_{13})_1(q s_{24})_1
(q s_{14})_1}\nonumber\\
&&\times
\varphi(1,1,2,1,1,1)\nonumber\\
&&+  {\zeta_4\zeta_4\over \zeta_1\zeta_2} \left(q {s_1s_2\over s_4s_4} \right) 
{(q t^{-1})_2(q t^{-1})_1(q t^{-1})_1
(t)_2(t)_1(t)_1  (q^2 \,s_{13})_1\over
(q)_1(q)_1(q s_{12})_1(q s_{23})_1
(q s_{34})_1(q s_{13})_1(q s_{24})_2
(q s_{14})_1}\nonumber\\
&&\times
\varphi(1,1,1,1,2,1)\nonumber\\
&&+  {\zeta_4\zeta_4\over \zeta_1\zeta_2} \left(-q{s_1s_2\over s_4s_4} \right) 
{(q t^{-1})_2(q t^{-1})_2(q t^{-1})_1
(t)_2(t)_2(t)_1  (q^{3} s_1s_2/s_3s_4)_1\over
(q)_1(q)_1(q s_{12})_1(q s_{23})_2
(q s_{34})_2(q s_{13})_1(q s_{24})_2
(q s_{14})_1}\nonumber\\
&&\times
\varphi(1,2,2,1,2,1),\nonumber
\end{eqnarray}
and
\begin{eqnarray}
&&y_{2,2,2} \\
&&= {\zeta_4\zeta_4\over \zeta_1\zeta_1} \left(q^2{s_1s_1\over s_4s_4} \right) 
{(q t^{-1})_2(q t^{-1})_2(t)_2(t)_2 \over
(q)_2(q s_{12})_2(q s_{24})_2(q s_{14})_2}
\varphi(2,0,0,0,2,2)\nonumber\\
&&+ {\zeta_4\zeta_4\over \zeta_1\zeta_1} \left(q^2{s_1s_1\over s_4s_4} \right) 
{(q t^{-1})_2(q t^{-1})_2(t)_2(t)_2 \over
(q)_2(q s_{13})_2(q s_{34})_2(q s_{14})_2}
\varphi(0,0,2,2,0,2)\nonumber\\
&&+ {\zeta_4\zeta_4\over \zeta_1\zeta_1} \left(q{s_1s_1\over s_4s_4} \right) 
{(q t^{-1})_2(q t^{-1})_2(q t^{-1})_2
(t)_2(t)_2(t)_2  (q^{3} \,s_1s_2/s_3s_4)_2\over
(q)_2(q s_{12})_2(q s_{23})_2
(q s_{34})_2(q s_{13})_2(q s_{24})_2
(q s_{14})_2} \nonumber\\
&&\times
\varphi(2,2,2,2,2,2) \nonumber\\
&&+ {\zeta_4\zeta_4\over \zeta_1\zeta_1} \left(-q{s_1s_1\over s_4s_4} \right) 
{(q t^{-1})_1(q t^{-1})_1(q t^{-1})_1
(t)_1(t)_1(t)_1  (q^2 \,s_2s_4/s_1s_3)_1\over
(q)_1(q s_{12})_1(q s_{23})_1
(q s_{34})_1(q s_{13})_1(q s_{24})_1
(q s_{14})_1} \nonumber\\
&&\times
\varphi(1,1,1,1,1,1) \nonumber\\
&&+ {\zeta_4\zeta_4\over \zeta_1\zeta_1} \left(q{s_1s_1\over s_4s_4} \right) 
{(q t^{-1})_2(q t^{-1})_1(q t^{-1})_1
(t)_2(t)_1(t)_1  (q^2 \,s_{23})_1\over
(q)_1(q)_1(q s_{12})_1(q s_{23})_1
(q s_{34})_1(q s_{13})_1(q s_{24})_1
(q s_{14})_2} \nonumber\\
&&\times
\varphi(1,1,1,1,1,2) \nonumber\\
&&+ {\zeta_4\zeta_4\over \zeta_1\zeta_1} \left(q{s_1s_1\over s_4s_4} \right) 
{(q t^{-1})_2(q t^{-1})_1(q t^{-1})_1
(t)_2(t)_1(t)_1  (q^2 \,s_{21})_1\over
(q)_1(q)_1(q s_{12})_1(q s_{23})_1
(q s_{34})_2(q s_{13})_1(q s_{24})_1
(q s_{14})_1} \nonumber\\
&&\times
\varphi(1,1,2,1,1,1) \nonumber\\
&&+ {\zeta_4\zeta_4\over \zeta_1\zeta_1} \left(q{s_1s_1\over s_4s_4} \right) 
{(q t^{-1})_2(q t^{-1})_1(q t^{-1})_1
(t)_2(t)_1(t)_1  (q^2 \,s_{43})_1\over
(q)_1(p^{1\over 2})_1(q s_{12})_2(p^{1\over 2}s_{23})_1
(q s_{34})_1(q s_{13})_1(q s_{24})_1
(q s_{14})_1}\nonumber\\
&&\times
\varphi(2,1,1,1,1,1)\nonumber\\
&&+  {\zeta_4\zeta_4\over \zeta_1\zeta_1} \left(-q{s_1s_1\over s_4s_4} \right) 
{(q t^{-1})_2(q t^{-1})_2(q t^{-1})_1
(t)_2(t)_2(t)_1  (q^{3} \,s_1s_2/s_3s_4)_1\over
(q)_1(q)_1(q s_{12})_1(q s_{23})_1
(q s_{34})_2(q s_{13})_2(q s_{24})_1
(q s_{14})_2}\nonumber\\
&&\times
\varphi(1,1,2,2,1,2)\nonumber\\
&&+  {\zeta_4\zeta_4\over \zeta_1\zeta_1} \left(-q{s_1s_1\over s_4s_4} \right) 
{(q t^{-1})_2(q t^{-1})_2(q t^{-1})_1
(t)_2(t)_2(t)_1  (q^{3} \,s_1s_2/s_3s_4)_1\over
(q)_1(q)_1(q s_{12})_2(q s_{23})_1
(q s_{34})_1(q s_{13})_1(q s_{24})_2
(q s_{14})_2}\nonumber\\
&&\times
\varphi(2,1,1,1,2,2)\nonumber\\
&&+  {\zeta_4\zeta_4\over \zeta_1\zeta_1} \left(-q {s_1s_1\over s_4s_4} \right) 
\nonumber\\
&&\times
{(q)_2
(q t^{-1})_2(q t^{-1})_2(q t^{-1})_1
(t)_2(t)_2(t)_1 \over
(q)_1(q)_1(q)_1(q)_1
(q s_{12})_2(q s_{23})_1
(q s_{34})_2(q s_{13})_1(q s_{24})_1
(q s_{14})_1}\nonumber\\
&&\times
\varphi(2,1,2,1,1,1).\nonumber
\end{eqnarray}
Here, we have used the notations $s_{ij}=s_i/s_j$,
$(a)_k=(a;q)_k$ and
\begin{eqnarray}
&&\varphi(k_{12},k_{23},k_{34},k_{13},k_{24},k_{14})\\
&=&
\prod_{1\leq i<j\leq 4}
(1-\zeta_j/\zeta_i)
\cdot {}_2\phi_1\left(
{  q^{k_{ij}+1} t^{-1}, q t^{-1}s_i/s_j ,
\atop
q^{{k_{ij}+1}}s_i/s_j  };q, t\zeta_j/\zeta_i
\right).\nonumber
\end{eqnarray} 
 \end{prop}
 
We may convince ourselves that the above expression gives us an efficient 
way to produce correct coefficients for the eigenfunction. 
Furthermore, 
even the following can be observed.
 \begin{con}
 The expression for the first eigenfunction given in Proposition \ref{n=4-eigen}
gives correct coefficients for the 
powers $ ({\zeta_2/ \zeta_1})^{i_1}
  ({\zeta_3/\zeta_2})^{i_2}
 ({\zeta_4/\zeta_3})^{i_3}$ satisfying the conditions:
  \begin{eqnarray}
&&0\leq  i_1<\infty ,\quad 0\leq  i_2\leq 2,\quad 0\leq  i_3<\infty,\\
{\rm or}&&
0\leq  i_1\leq 2 ,\quad 0\leq  i_2<\infty,\quad 0\leq  i_3\leq 2.
 \end{eqnarray}
 \end{con}
 
The author has extended the formula given in 
Proposition \ref{n=4-eigen} for a little higher powers in $\zeta_i$'s, and observed that
there still exists  a nice hypergeometric-like structure. At this moment, 
however, it is not clear how to organize the general terms, 
and we omit writing.
\bigskip

In the first several terms of the
eigenfunctions given in Proposition \ref{n=4-eigen},
we do not see any dependence on the parameter $\alpha$. 
Thus, Conjecture \ref{Con1} is very likely true for $n=4$.

\section{Concluding Remarks}
Let us summarize the results obtained in this paper.
In Section 1, we have introduced the integral transformation $I(\alpha)$, and
it was conjectured that $I(\alpha)$ generates a commutative family of 
operators acting on the space of formal series ${\cal F}_n$ (Conjecture \ref{Con1}).
In Section 2, a proof of the commutativity for the case $n=2$ was given
(Theorem \ref{Thm2}). In Section 3, the existence of the eigenfunctions
for general $n$ was studied, when all the parameters are generic 
(Proposition \ref{existing}).
Explicit formulas for the eigenfunctions were examined for 
$n=3,4$ (Conjecture \ref{eigen_for_n=3} and Proposition \ref{n=4-eigen}).
Since these eigenfunctions 
do not depend on the spectral parameter $\alpha$ of the
integral transformation $I(\alpha)$, it is expected that 
Conjecture \ref{Con1} is true for $n=3,4$.
\bigskip

Finally, let us make some comments to show 
the directions toward the next papers \cite{S2,S3}. 
When we specialize the parameters as
$s_1=s_2=\cdots=s_n$,
there exist several phenomena
which are interesting both from the 
hypergeometric series and the lattice model viewpoints.
In \cite{S2}, a class of hypergeometric-type functions will be introduced,
which is characterized by 
a certain covariant transformation property 
with respect to the action of $I(\alpha)$ together with
an initial condition given by an infinite product expression.
To be more precise, by setting $s_1=s_2=\cdots=s_n=1$,
we will introduce a series 
$F(\alpha)=
F(\zeta_1,\cdots,\zeta_n;\alpha,q,t)\in {\cal F}_n$ 
characterized by the conditions
\begin{eqnarray}
{\rm (I)}&&
I(\alpha q^{-1}t ) \cdot F(\alpha)=F(\alpha q^{-{1}}t),\\
{\rm (II)}&&
F(t^{1/2})
=\prod_{1\leq i<j\leq n}(1-\zeta_j/\zeta_i)
{(q t^{-1/2}\zeta_j/\zeta_i;q)_\infty \over 
(t^{1/2}\zeta_j/\zeta_i;q)_\infty }.
\end{eqnarray}
This function $F(\alpha)$ will be called
`quasi-eigenfunction' for short.
We are interested in the function $F(\alpha)$ because of the
following observations:
\begin{enumerate}
\item The explicit formula of $F(\alpha)$ for $n=2$ can be easily derived \cite{S2}.
Further, 
we are able to have a conjectural expression of $F(\alpha)$ for $n=3$:
\begin{eqnarray}
&&F(\zeta_1,\zeta_2,\zeta_3;\alpha,q,t)\\
&=&
\sum_{k=0}^\infty
{
(\alpha^{-2}t,q t^{-1},q t^{-1};q)_k \over 
(q,\alpha^{-1}q,\alpha^{-1}q;q)_k}
(q \zeta_3/\zeta_1)^k
{}_2\phi_1 \left({\alpha^{-1}, q^{-{k}} \atop \alpha q^{-k+1}};
q, \alpha t\right)
\nonumber\\
&&\times
\prod_{1\leq i<j\leq 3}
(1-\zeta_j/\zeta_i)\;
{}_2\phi_1\left(
{  q^{k+1} t^{-1}, \alpha q t^{-1}
\atop
\alpha^{-1}q^{{k+1}} };q, \alpha^{-1}t\zeta_j/\zeta_i
\right).\nonumber
\end{eqnarray}
\item
We may find a variety of infinite product expressions for $F(\alpha)$
when we suitably specialize the parameter $\alpha$.
Simplest examples among these are:
\begin{eqnarray}
&&
F(- t^{1/2})
=\prod_{1\leq i<j\leq n}(1-\zeta_j/\zeta_i)
{(- q t^{-1/2}\zeta_j/\zeta_i;q)_\infty \over 
(- t^{1/2}\zeta_j/\zeta_i;q)_\infty }, \label{degen1}\\
&&
F( t )
=\prod_{1\leq i<j\leq n\atop {\rm step} 2}(1-\zeta_j/\zeta_i)
{( q t^{-1 }\zeta_j/\zeta_i;q)_\infty \over 
( t \zeta_j/\zeta_i;q)_\infty }. \label{degen2}
\end{eqnarray}
Here, we have used the notation
 $\prod_{1\leq i<j\leq n\atop {\rm step} 2} f_{ij}=
f_{13}f_{15}\cdots f_{24}f_{26}\cdots$.

\item
The highest-to-highest matrix elements of the eight-vertex vertex operators 
can be realized as 
\begin{eqnarray}
&&\langle \Phi(\zeta_1)\Phi(\zeta_2)\cdots \Phi(\zeta_n)\rangle \\
&=&
\prod_{1\leq i<j\leq n}
{\xi(\zeta_j^2/\zeta_i^2;p,q) \over 1-\zeta_j/\zeta_i} \cdot
F(\zeta_1,\cdots,\zeta_n;-1,p^{1/2},q),\nonumber
\end{eqnarray} 
where
\begin{eqnarray}
&&
\xi(z;p,q)={(q^2z;p,q^4)_\infty (pq^2z;p,q^4)_\infty \over 
(q^4z;p,q^4)_\infty (pz;p,q^4)_\infty }.
\end{eqnarray}
\end{enumerate}
One finds that the integral formulas of the matrix elements 
of the vertex operators for the cases $p^{1/2}=q^{3/2}$ and 
$p^{1/2}=-q^2$ given in \cite{S} can be easily recovered from the 
above properties of $F(\alpha)$. Note, however, the case $p^{1/2}=q^3$ corresponds to 
another type of product formula for $F(\alpha)$ which is not given here.

In the continuation of the present work \cite{S2,S3}, basic properties 
of the quasi-eigenfunction $F(\alpha)$ will be discussed.
\bigskip

\noindent
{\it Acknowledgment.}~~~
This work is supported by the  Grant-in-Aid for Scientific Research 
(C) 16540183.

\end{document}